\documentclass[12pt]{amsart}
\usepackage{amssymb}
\usepackage{amsmath}
\usepackage{amscd}
\usepackage[all]{xy}

\usepackage{color}
\date{\today}

\def\sL{{\mathcal L}}
\def\beq{\begin{equation} \label}

\def\sO{{\mathcal O}}

\newtheorem{theo}{Theorem}[section]
\newtheorem{prop}[theo]{Proposition}
\newtheorem{lem}[theo]{Lemma}
\newtheorem{cor}[theo]{Corollary}
\newtheorem{defi}[theo]{Definition}
\newtheorem{conj}[theo]{Conjecture}
\newtheorem{rem}[theo]{Remark}
\newtheorem{que}[theo]{Question}

\newcommand{\bthe}{\begin{theo}}
\newcommand{\ble}{\begin{lem}}
\newcommand{\bpr}{\begin{prop}}
\newcommand{\bco}{\begin{cor}}
\newcommand{\bde}{\begin{defi}}
\newcommand{\ethe}{\end{theo}}
\newcommand{\ele}{\end{lem}}
\newcommand{\epr}{\end{prop}}
\newcommand{\eco}{\end{cor}}
\newcommand{\ede}{\end{defi}}
\newcommand{\bconj}{\begin{conj}}
\newcommand{\econj}{\end{conj}}
\newcommand{\brem}{\begin{rem}}
\newcommand{\erem}{\end{rem}}
\newcommand{\bque}{\begin{que}}
\newcommand{\eque}{\end{que}}

\def \to {{\rightarrow}}

\def \F {{\mathbb F}}

\def \G {{\mathbb G}}

\def\k{{\overline k}}

\begin{document}

\baselineskip=16pt

\title{Potential density  for some families of homogeneous spaces}

\author{J.-L. Colliot-Th\'el\`ene}
\author{J. N. Iyer}
\address{CNRS, UMR 8628, Math\'ematiques, B\^atiment 425, Universit\'e Paris-Sud,
F-91405 Orsay, France}
\email{jlct@math.u-psud.fr}
\address{The Institute of Mathematical Sciences, CIT
Campus, Taramani, Chennai 600113, India}
\address{Department of Mathematics and Statistics, University of Hyderabad, Gachibowli, Central University P O, Hyderabad-500046, India}
\email{jniyer@imsc.res.in}
\date{February 21st, 2011}

\begin{abstract}
For a smooth, projective  family of homogeneous varieties defined over a number field, 
we show that if potential density holds for the
rational points of the base, then it also holds for the total space. A conjecture of  Campana and Peternell, known in dimension
at most 4 and for certain higher dimensional cases,
 would then imply potential density for the rational points of smooth projective varieties over number fields
whose tangent  bundle is nef.
\end{abstract}
\footnotetext{Mathematics Classification Number: 14C25, 14D05, 14D20, 14D21 }
\maketitle

\section*{Introduction}

Let $k$ be  a number field.   A geometrically integral variety $X$ over the field $k$ satisfies potential density if there exists a
finite field extension $K/k$ such that the set $X(K)$ of rational points of $X$ is Zariski dense in $X_{K}=X \otimes_{k}K$.
One hopes that this property only depends on the geometry of the variety $X$ over an algebraically closed field
containing $k$, for instance over the complex numbers.
It has been known for some time that Abelian varieties satisfy potential density  (see \cite[Prop. 4.2]{Hassett}).
It  is an open problem whether potential density holds for rationally connected varieties, in particular for Fano varieties.

   For an overview of problems and results regarding potential density over number field, as of 2003, including work
   of Bogomolov, Hasssett, Tschinkel,
   we refer  the reader to
the survey \cite{Hassett} by B. Hassett.    Among the significant later results, let us mention the paper by E. Amerik and C. Voisin \cite{AmerikVoisin}.

  According to  the Hartshorne--Frenkel conjecture, proved by S. Mori, a smooth, projective, complex variety whose
  tangent bundle is ample, is isomorphic to projective space. Over an arbitrary field $k$ of characteristic zero, this implies
  that a smooth, projective, geometrically integral $k$-variety $X$ whose tangent bundle is ample is a Severi-Brauer variety. After a finite extension  $K/k$ of  the ground field, this variety acquires a rational point and then it $k$-isomorphic to projective space over $K$, hence the set $X(K)$ is  Zariski dense in $X_{K}$. This argument of course has nothing to do with number fields.
  
  One  may wonder whether potential density holds more generally for a smooth, projective, geometrically integral $k$-variety $X$
  whose tangent bundle is nef. Such varieties have been studied in particular  by Campana,  Demailly, Peternell, Schneider.
  In this note we give a detailed proof of  a  stability property for potential density   (Theorem \ref{maintheorem} and Corollary \ref{applicab}). The result should be more or less obvious  to experts. Combined with a conjecture of Campana and Peternell,
  it  predicts potential density for varieties with nef tangent bundle.

\section{Known results on homogeneous spaces of linear algebraic groups}

The following theorem gathers results of T. A. Springer, J.-C. Douai and M. Borovoi (\cite{Bv}).

\bthe\label{generalhomogeneous}
Let $k$ be a field of characteristic zero and let $\k$ be an algebraic closure of $k$.
Let $G/k$ be a semisimple simply connected group.
Let $X/k$ be a homogeneous space of $G$. Assume that
a geometric stabilizer $\overline H$   is connected.

(a) The homogeneous space structure on $X$  defines a
$k$-kernel $\sL:=(\overline H,\kappa)$, and  a class $\eta(X)$
in the cohomology set $H^2(k,\sL)$. This class is neutral
if and only if  there exists  a principal homogeneous space $E$ under $G$
and a $G$-equivariant map $E \to X$.

(b) Let  ${\overline H}^{tor}$ be the maximal toric quotient of $\overline H$.
The $k$-kernel $\sL$ induces a $k$-kernel $({\overline H}^{tor},\kappa^{tor})$.
To the latter is associated a natural $k$-torus $T$. There is an induced map
of sets
$$H^2(k,\sL) \to H^2(k,T).$$
Let $\eta^{tor}_{X} \in H^2(k,T)$ denote the image of $\eta(X)$.

(c) If $L/k$ is a finite field extension such that $X(L) \neq \emptyset$
then $$[L:k].\eta^{tor}_{X}=0 \in H^2(k,T).$$

(d) If $X/k$ is projective, then $\overline H$ is connected, and the
associated torus $T$ is a quasitrivial torus.
\ethe

Proof.  For (a), (b), (c), see \cite{Bv} and the review in  \cite[\S 5, p. 333--335]{CTGP}.
For (d), see \cite[Lemma 5.6]{CTGP}.

The following theorem combines results of Kneser, Bruhat-Tits (for principal homogeneous spaces of
semisimple simply connected groups) and Springer, Douai, Borovoi.
\bthe\label{homogeneouspadic}
Let $k$ be a $p$-adic field. In the situation of Theorem~\ref{generalhomogeneous},
the class $\eta(X)$ is neutral if and only if $\eta^{tor}_{X}=0 \in H^2(k,T)$.
In that case, $X$ has a $k$-point.
\ethe

Proof. \cite[Thm. 5.5]{Bv} and   \cite[Prop. 5.4]{CTGP}.

\bpr\label{stabparabo}
Let $k$ be a field of characteristic zero and let $\k$ be an algebraic closure of $k$.
 Let $X/k$ be a smooth, projective, geometrically connected variety. 

(a) If $X\times_{k}\k$ is a homogeneous space of a linear algebraic group, then
there exists a semisimple simply connected group $G$ over $k$
such that $X$ is a homogeneous space of $G$. 

(b)The geometric
stabilizers of this action are parabolic groups, in particular they
are connected.

(c) If $X(k) \neq \emptyset$, then $X$ is $k$-birational to projective space.
\epr

Proof. Statement (a) is a special case of the following theorem of Demazure.
The idea here is to consider the neutral component $G=Aut^0_{X/k}$ of the automorphism
group of $X$ over $k$, which is an adjoint group, and then to take the semisimple cover
of that group.
For (b), see \cite[IV.11.6]{borel}. For (c), see \cite[IV.14.21 and V.20.5] {borel}.
\bigskip

\bthe\label{homogeneousfamily} 
Let $k$ be a field of characteristic zero and let $\k$ be an algebraic closure of $k$.
Let $p : X \to Y$ be a  smooth, proper $k$-morphism of smooth, geometrically connected $k$-varieties.

(a) (Demazure) If each geometric fibre of $p$ is a homogeneous space of a connected  linear algebraic group then
the group $G=Aut^0_{X/Y}$  is a semisimple group over $Y$ and $X \to Y$ is
a homogeneous space of $G$. The fibres of $G$ are adjoint groups. There also exists
a semisimple group $G^{sc}$ over $Y$, whose fibres are simply connected semisimple groups, 
such that $X$ is a homogeneous space of $G^{sc}$.

(b) In the above situation, there exists a  finite Zariski open cover $\{U_{i} \}_{i \in I}$ of $Y$
and quasifinite, surjective, \'etale maps $V_{i} \to U_{i}$ which factorize as
$V_{i} \to X\times_{Y}U_{i}  \to U_{i}.$ 

(c) There exists an integer $d$ such that for any point $M \in Y$ 
the torus $T_{M}$ over the residue field $k(M)$ associated
to the homogeneous space $X\times_{Y}M$ (see Theorem \ref{generalhomogeneous} (b))
has rank $d$.

(d) There exists an integer $N>0$ such that
for any field $L$ containing $k$ and any $L$-point $P \in Y(L)$,
the class $\eta^{tor}_{X_{P}} \in H^2(L,T_{P})$ is $N$-torsion.
Here $T_{P}$ denotes the $k$-torus associated to 
the $k$-variety  $X_{P}$ (fibre of $X \to Y$ at $P$)  viewed as a homogeneous space
of $G^{sc}_{P}$.
\ethe

Proof. Statement (a) is   \cite[Prop. 4]{Demazure}) of Demazure.
Statement (b)  is a general fact for a smooth, surjective morphism $X \to Y$.
For any point $P\in Y(k)$ there exists an $i$ with $P \in U_{i}$ 
and a closed point $M \in V_{i}$ mapping to $P$. 
Let $k(M)$ be the residue field of $M$.
Since the set $I$ is finite  and for each $i$ the degrees
of the fibres of $V_{i} \to U_{i}$ are bounded, 
there exists a fixed integer $N>0$, independent of $P$,  such that
the degree of the field extension
$k(M)/k$ divides   $N$.
We now use Theorem \ref{generalhomogeneous}, which applies to the present situation
in view of
Proposition \ref{stabparabo}. The class $\eta \in H^1(k,T_{P})$ vanishes in
$H^1(k(M),T_{P})$. Hence its corestriction $[k(M):k].\eta$ vanishes in $H^1(k,T_{P})$.
So does $N.\eta$.

\bthe (Harder) \label{harder}
Let $k$ be  a number field. Let $X/k$ be a smooth projective
homogeneous variety under the action of a connected  linear algebraic group.
Then the Hasse principle holds for $X$ : if $X$ has points in all
completions of $k$, then it has a point in $k$. 
\ethe

Proof. In \cite{Harder}, Harder reduces the local-global statement
to the Hasse principle for principal homogeneous spaces of 
semisimple simply connected groups. In this set-up, the local-global 
principle is due to Eichler, Kneser, Harder, and Chernousov.

\section{The theorem}\label{proofoftheorem}

To prove the main theorem, we shall use two further results. The first one is
a special case of  a standard result in the study of the Hasse principle.
  
\bthe\label{almostonto}
 Let $k$ be a number field. Let $p : X \to Y$ be a smooth, projective morphism
of projective, geometrically integral $k$-varieties.
Assume that the fibres of $p$ are
homogeneous spaces of connected linear algebraic groups. 
Then there exists a finite set $S$ of places of $k$ such that
for any finite field extension $L/k$ and any place $w$ of $L$
not lying above a place in $S$ the induced map
$X(L_{w}) \to Y(L_{w})$ is onto.
\ethe

Proof.  By Theorem \ref{homogeneousfamily}, there exists a semisimple group $G$ over $Y$
such that $X$ is a homogeneous space of $G$.
By a standard limit argument, which is easy in the present, projective context,
 (for a more general set up, see EGA IV 8), the whole situation may be spread out over
an open set  $B$ of the spectrum of the  ring of integers of $k$. Let ${\bf X} \to {\bf Y}$
and ${\bf G}/ {\bf Y}$ denote the corresponding objects.
Let $v$ be a place in $B$. Let $O_{v} \subset k_{v}$ denote the
ring of integers in the completion $k_{v}$, and let $\F_{v}$ denote
the residue field. Let $P_{v} \in Y(k_{v})$.
Since ${\bf Y}/B$ is proper,
we have ${\bf Y}(O_{v})=Y(k_{v})$, the point $P_{v}$ may be viewed
as a point  ${\bf P}_{v} \in {\bf Y}(O_{v})$. By restriction to ${\bf P}_{v}$
one gets a homogeneous space of the $O_{v}$-semisimple group
${\bf G}\times_{Y}{\bf P}_{v}$.  One then considers the reduction of all this
over the finite field $\F_{v}$.
Any homogeneous space of a connected linear algebraic group
over a finite field  has a rational point (Lang, Springer, see Serre \cite[Chap. III, \S 2]{serre}). 
By Hensel's lemma one then
lifts such a point to an $O_{v}$-point of ${\bf X}\times_{\bf Y}{\bf P}_{v}$.
Such a point defines a $k_{v}$-point of $X$ whose image is $P_{v} \in Y(k_{v})$.
Thus $X(k_{v}) \to Y(k_{v})$ is onto. The same argument works over any
finite field extension $L$ of $k$, with the inverse image of $B$ in the spectrum
of the ring  of integers of $L$.

\ble\label{padiclemma}
Let $k$ be a $p$-adic field. Let $T$ be a quasisplit torus
of dimension $d$. Let $N>0$ be an integer. 
If $L/k$ is a field extension whose degree is divisible by $N.d!$
then the restriction map on $N$-torsion classes
$$ H^2(k,T)[N] \to H^2(L,T)[N]$$
is zero.
\ele

Proof. We immediately reduce to the case $T=R_{K/k}\G_{m}$, where $K/k$
is a field extension of degree $r \leq d$. By a lemma of Faddeev and Shapiro (\cite[Chap. I, \S 2.5]{serre}), the restriction map
$H^2(k,T) \to H^2(L,T)$ then reads
$Br(K) \to \oplus_{i} Br(L_{i}),$
where $L\otimes_{k}K = \prod_{i} L_{i}$ is the decomposition into  a finite product of
fields. We have the embeddings $k \subset K \subset L_{i}$ and $k \subset L \subset L_{i}$.
By assumption, $N.d!$ divides $ [L:k]$, which  divides  $[L_{i}:k]= [K:k][L_{i}:K]=r[L_{i}:K]$.  It follows that $N$ 
divides $[L_{i}:K]$. But the map of Brauer groups of local fields
$Br(K) \to  Br(L_{i})$ reads as multiplication by $[L_{i}:K]$ on ${\bf Q}/{\bf Z}$ (\cite[Chap. XIII, \S 3, Prop. 7 p.~201]{corpslocaux}.
Hence on $N$-torsion it is zero.

\bthe\label{maintheorem}
Let $k$ be a number field. Let $p : X \to Y$ be a smooth, proper morphism
of geometrically integral varieties. Assume that the geometric fibres of $p$ are
homogeneous spaces of  connected linear algebraic groups. Then there exists a finite field extension $L/k$
such that $Y(k)  \subset Y(L)$ lies in the image of $ X(L) \to Y(L)$.
If $Y(k)$ is Zariski dense in $Y$, then  for $L$ as above, $X(L)$ is Zariski dense in $X_L$.
\ethe

Proof.  By Theorem \ref{almostonto}, there exists a finite set $S$ of places of $k$,
which we assume to contain all archimedean places, such that 
for any finite field extension $L/k$ and any place $w$ of $L$
not lying above a place in $S$ the induced map
$X(L_{w}) \to Y(L_{w})$ is onto. By Theorem
\ref{homogeneousfamily}, there exists an integer $d>0$
 and an integer $N>0$, which we may choose even,  such that
for any field $L$ containing $k$, and any point $M \in Y(L)$,
the torus $T_{M}$ over $L$ associated to the homogeneous space $X_{M}$ defined by the fibre at $M$
is a quasitrivial torus  over field $L$, of dimension $d$, and the class $\eta^{tor}(X_{M}) \in H^2(L,T_{M})$ is
annihilated by $N$. For each finite place $v \in S$ let us pick a field extension $F^{v}/k_{v}$ of degree
$N.d!$. For each archimedean place $v$ of $k$ let $F^v/k_{v}$ be a separable extension of $k_{v}$
of degree $N.d!$, hence even, which breaks up as the product of copies of the complex field.
 By weak approximation for the field $k$ and Krasner's lemma \cite[Chap. II, \S 2, Exercice 2, p.~40 ]{corpslocaux},
  there exists a field extension
 $L/k$ of degree $N.d!$ such that for each $v \in S$, there is an isomorphism $L\otimes_{k}k_{v} \simeq F^v$.
 In particular, for each finite place $v$ of $k$ in $S$, there is just one place $w$ of $L$ above $v$.
 
 Let now $P \in Y(k)$ be an arbitrary point, let $T=T_{P}$ be the $k$-torus of dimension $d$  associated to the
 homogeneous space $X_{P}$
 and let $\eta= \eta^{tor}(T) \in H^2(k,T)$ be
 the associated class.  
 This  class  is annihilated by $N$. At any place $v$ of $k$ not in $S$,
 the fibre $X_{P}$ has a $k_{v}$-point, hence $\eta_{v}=0 \in H^2(k_{v},T)$.
 If $w$ is a place of $L$ over a place of $S$, Lemma \ref{padiclemma}
 and the choice of the extension $L/k$ imply that the image of $\eta$ in $H^2(L_{w},T)$
 vanishes. Thus $\eta_{L} \in H^2(L,T)$ vanishes over each completion of $L$.
By theorem \ref{homogeneouspadic} this implies that $X_{P}\otimes_{k}L$ has points in all completions of $L$.
By Theorem \ref{harder} this implies that the $L$-variety $X_{P}\otimes_{k}L$ has an $L$-point,
and then that $X_{P}\otimes_{k}L$  is $L$-birational to projective space over $L$, in particular 
$L$-points are Zariski dense on $X_{L}$. This completes the proof of the theorem.

\medskip

\brem
{\rm  In the more general context of integral points,
a special case of the above theorem (family of Severi--Brauer varieties)
was remarked some time ago  by the first named author  \cite[Thm. 2.8]{HT}.
One could certainly also write down an integral 
points version of Theorem~\ref{maintheorem}.}
\erem

\bco\label{applicab}
Let $k$ be a number field. Let $A$ be abelian variety over $k$.
Let $p : X \to A$ be a smooth, proper morphism
of geometrically integral varieties. Assume that the geometric fibres of $p$ are
homogeneous spaces of  connected linear algebraic groups.
Then there exists a finite field extension $K/k$ such that $X(K)$ is Zariski dense in $X_{K}$.
\eco
 
 Proof. Since potential density holds for abelian varieties (\cite[Prop.~4.2]{Hassett}), this
 is an immediate consequence of  Theorem  \ref{maintheorem}.

\section{Varieties with nef tangent bundles : the conjecture of Campana and Peternell}

In this section we discuss
potential density of rational points for smooth, projective, geometrically integral varieties over number field,
under the assumption that their tangent bundle is
  numerically effective (nef).
  By definition, this  means that the line bundle $\sL:=\sO_{PT(X)}(1)$ on the projectivized tangent bundle $PT(X)$, is numerically effective, i.e. $\sL.C\geq 0$, for any curve $C$ on $PT(X)$.

Recall that a smooth, projective variety $X$ is a Fano variety if the anticanonical line bundle $-K_{X}$
is ample.
The following theorem was conjectured by Campana and Peternell and proved by them in
dimension at most 3  \cite[Theorem, p.169]{Ca-Pe}.

\begin{theo}(Demailly--Peternell--Schneider) \label{dps}  \cite[Main Theorem, p.~296]{DPS} 
Let $k$ be an algebraically closed field of characteristic zero. Let $X$ be 
a smooth, projective, connected variety with nef tangent bundle. Then there exists
a finite \'etale connected cover   $X' \to X$ such that for any $k$-point of $X'$
the associated Albanese map $X' \to A$ to the Albanese variety of $X'$ (which is an 
abelian variety) is a smooth, projective morphism whose fibres are 
Fano varieties with nef tangent bundles.
\end{theo}

Campana and Peternell  put forward the following conjecture.

\bconj\label{conjCaPe}
  \cite[Conjecture 11.1, p.~185]{Ca-Pe} Over an algebraically closed field of characteristic zero,
a Fano variety with nef tangent bundle is a  projective homogeneous variety of a linear
algebraic group, i.e. it is of the shape $G/P$ for $G$ a connected linear algebraic group
and $P$ a parabolic subgroup.
\econj

  A variant is formulated by J-M. Hwang   \cite[Conjecture 4.1, p.~622]{Hwang}:
this should be the case as soon as  all rational curves on  $X$ are free.

The Campana-Peternell conjecture \ref{conjCaPe} was proved by Campana and Peternell  in dimension up to  3
and by J.-M. Hwang in dimension 4. It has also been proved  
 for higher dimensional Fano, when the Betti numbers satisfy $b_2=b_4=1$, 
 and the variety of minimal rational tangents at a general point is one-dimensional \cite[Main Theorem, p.~2641]{Mok}, \cite[Theorem 4.3, p.~623]{Hwang}. 
See \cite[section 4]{Hwang} for a discussion and references. See also a related recent work \cite{Biswas}. 
In  these various  cases, the following theorem therefore applies.

\bthe\label{condtheo}
Suppose $X$ is a smooth projective variety with a nef  tangent bundle, defined over a number field.
 Under Conjecture  \ref{conjCaPe} on Fano varieties, potential density holds for $X$.
\ethe

Proof. Combine Theorem  \ref{dps} (which descends from an
algebraic closure of $k$ to some finite extension of $k$), 
 Conjecture \ref{conjCaPe} and Corollary \ref{applicab}.

\bque
Let $Y$ be an abelian variety  over an algebraically closed field $k$. Let $X \to Y$
be a smooth, projective family of homogenous spaces of connected linear algebraic groups.
Does there exist a finite \'etale map  $Z \to Y$ with $Z$ 
connected such that $X\times_{Y}Z \to Z$ admits a rational section?
\eque

Since potential density is inherited by  finite \'etale covers (Chevalley--Weil, cf. \cite[Prop. 3.4]{Hassett}),
an affirmative answer to the question would lead to 
 an alternate proof of Theorem \ref{maintheorem}.

Over an algebraically closed field, a connected, finite \'etale cover of an abelian variety
may be given the structure of an abelian variety.
If  
the above question had an affirmative answer,
this    would give  an alternate, less arithmetic  proof for  
Corollary \ref{applicab} and therefore for Theorem \ref{condtheo}.

In the special case where $X \to Y$ is a Severi-Brauer scheme,
the answer to the above question is in
the affirmative (see the proof of  \cite[Lemma 7.4 (1)]{Ca-Pe}).

\bigskip

  {\Small Acknowledgement: This work was discussed at the Universit\'e  Paris-Sud, Orsay, during Nov-14-Dec 14 2010. The second named author thanks Laurent Clozel for the kind invitation and the D\'epartement de Math\'ematiques, for their hospitality and support from the ARCUS programme.
}


\begin{thebibliography}{AAAAA}

\bibitem[Am-Vo]{AmerikVoisin} 
E. Amerik and C. Voisin,   {\em   Potential density of rational points on the variety  of lines of a cubic fourfold}, Duke Math. J., \textbf{145}, No. 2,  379-408,  (2008).

\bibitem[Bi-Br]{Biswas} I. Biswas and U. Bruzzo, {\em Holomorphic Cartan geometry on manifolds with numerically effective tangent bundle}, arXiv math.AG.1101.4192.

\bibitem[Bo]{borel} A. Borel, Linear algebraic groups, Second enlarged edition, GTM {\bf 126} Springer-Verlag (1991).

\bibitem[Bv]{Bv} 
M. Borovoi,  {\em Abelianization of the second nonabelian Galois cohomology}, Duke Math. J., \textbf{72}, (1993), 217-239.


\bibitem[Ca-Pe]{Ca-Pe} F. Campana and T. Peternell,  {\em Projective manifolds whose tangent bundles are numerically effective}, Math. Ann. \textbf{289} (1991), 169-187.


\bibitem[CTGP]{CTGP} J.-L. Colliot-Th\'el\`ene, P.  Gille and R. Parimala,  {\em Arithmetic of linear algebraic groups over 2-dimensional geometric fields}, Duke Math. J. \textbf{ 121}, No.2, 2004, 285-341.

\bibitem[DPS]{DPS} J.-P. Demailly, T. Peternell and M. Schneider,  {\em Compact complex manifolds with numerically effective tangent bundles}, Journal of Algebraic Geometry \textbf{3} (1994), 295-345.

\bibitem[De]{Demazure} M. Demazure,  {\em Automorphismes et d\'eformations des vari\'et\'es de Borel}, Invent. math. \textbf{39}, 179-186.

\bibitem[Ha]{Harder}  G. Harder,   {\em Bericht  \"uber neuere Resultate der Galoiskohomologie halb\-einfacher Gruppen}, Jber. Deutsch. Math.-Verein. \textbf{70} (1967/1968), 182--216.

\bibitem[Has]{Hassett} B. Hassett, Potential density of rational points on algebraic varieties, in Higher dimensional varieties
and rational points, ed. J. B\"or\"ocsky Jr., J. Koll\'ar, T. Szamuely, Bolyai Society Mathematical Studies \textbf{12} (2003), Springer. 
Budapest.

\bibitem[HT]{HT} B. Hassett and Yu. Tschinkel,  {\em Density of integral points on algebraic varieties}, 
in Rational points on algebraic varieties (E. Peyre, Yu. Tschinkel ed.),
p. 169--197,
Progress in Math. \textbf{199} (2001) Birkh\"auser Verlag.

\bibitem[Hw]{Hwang} J.-M. Hwang,   {\em Rigidity of rational homogeneous spaces}, Proceedings of ICM 2006, Madrid, Volume \textbf{2}, 613-626.

\bibitem[Mk]{Mok} N. Mok, {\em On Fano manifolds with nef tangent bundles admitting $1$-dimensional varieties of minimal rational tangents},
Trans. Amer.Math.Soc, \textbf{354}, (2002), 2639-2658.

\bibitem[S1]{corpslocaux} J.-P. Serre, Corps locaux, Publications de l'Institut Math\'ematique de l'Universit\'e de Nancago, VIII, Actualit\'es scientifiques et industrielles 1296, Hermann, Paris, 1968.

\bibitem[S2]{serre} J.-P. Serre,  Cohomologie galoisienne, cinqui\`eme \'edition, r\'evis\'ee et compl\'et\'ee,   LNM {\bf 5} (1973, 1994),
Springer Verlag.


\end{thebibliography}
\end{document}